\documentclass{jnmp01b}

\usepackage{amsmath}

\numberwithin{equation}{section}

\newtheorem{theorem}{Theorem}
\theoremstyle{definition}
\newtheorem{remark}{Remark}

\def\Real{\mathbb{R}}
\def\bb{\begin{equation}}
\def\ee{\end{equation}}
\def\mod{\hbox{mod}}
\def\pt{\partial}
\def\const{\hbox{const}}
\def\th{\hbox{tanh}}
\def\ch{\hbox{cosh}}
\def\a{\alpha}
\def\G{\Gamma}
\def\g{\gamma}
\def\l{\lambda}
\def\t{\tau}
\def\ve{\varepsilon}
\def\s{\sigma}

\makeatletter
\DeclareRobustCommand{\primfrac}[1]{%
  \PackageWarning{amsmath}{%
Foreign command \@backslashchar#1; %
\protect\frac\space or \protect\genfrac\space should be used instead%
  }
  \global\@xp\let\csname#1\@xp\endcsname\csname @@#1\endcsname
  \csname#1\endcsname
}
\makeatother

\begin{document}

\renewcommand{\evenhead}{O.M. Kiselev}
\renewcommand{\oddhead}{Dromion Perturbation  for the Davey-Stewartson-1
Equations}


\thispagestyle{empty}

\begin{flushleft}
\footnotesize \sf
Journal of Nonlinear Mathematical Physics \qquad 2000, V.7, N~4,
\pageref{firstpage}--\pageref{lastpage}.
\hfill {\sc Letter}
\end{flushleft}

\vspace{-5mm}

\copyrightnote{2000}{O.M. Kiselev}

\Name{Dromion Perturbation  for the Davey-Stewartson-1
Equations}

\label{firstpage}

\Author{O.M. KISELEV}

\Adress{Institute of Mathematics, Ufa Sci Centre of Russian Acad. of Sci,
112, \\
Chernyshevsky str., Ufa, 450000, Russia\\
E-mail: ok@imat.rb.ru}

\Date{Received December 23, 1999; Revised April 4, 2000;
Accepted June 13, 2000}

\begin{abstract}
\noindent
The perturbation of the dromion of the Davey-Stewartson-1 equation is
studied over  the large time.
\end{abstract}


\section{Introduction}

In this work we construct the asymp\-to\-tic so\-lu\-ti\-on of the
perturbed Da\-vey-Ste\-wart\-son-1 equations (DS-1):
\begin{equation}
\begin{split}
i\pt_t Q+{1\over2}(\pt_\xi^2+\pt_\eta^2)Q+(G_1+G_2)Q=\ve
iF,
\\
\pt_\xi G_1=-{\s\over2}\pt_\eta |Q|^2,\quad
\pt_\eta G_2=-{\s\over2}\pt_\xi |Q|^2.
\end{split}
\label{ds1}
\end{equation}
Here $\ve$ is small positive parameter, $\s=\pm1$ correspond to so
called focusing or defocusing DS-1 equations.

The equations (\ref{ds1}) describe the interaction between long and
short waves on the liquid surface in case the capillary effects and
potential flow are taken into account \cite{D-S,D-R}. The existence
theorems for the solutions of this equations in the various
functional classes are well-known \cite{Gh-S,Sh}. The inverse
scattering trans\-form method  for the DS equations was formulated in
\cite{Nizhnik}--\cite{F-Sung}. This method allows  to construct
solitons \cite{B-L-M-P} and to study global properties of the generic
solutions, for instance, the asymptotic behaviour over large times
\cite{OK1,OK2}. The asymptotics for nonintegrable DS equations were
studied in \cite{Nakao}.

The perturbations of the equations (\ref{ds1}) arise due to a small
irregularity of bottom  or by taking into account the next
corrections in more realistic models for the liquid surface
considered in \cite{D-S, D-R}. For the first case the perturbation
takes the form: $F\equiv AQ$. Here $A$ is real constant and its sign
corresponds to decreasing or increasing depth with respect to
spatial variable $\xi$.

We start with the solution of DS-1 constructed in \cite{B-L-M-P}:
\begin{equation}
q(\xi,\eta,t;\rho)=
{\rho\l\mu\exp(it(\l^2+\mu^2))\over2\cosh(\mu
\xi)\cosh(\l\eta)(1-\s|\rho|^2{\mu\l\over16}(1+\th(\l\eta))
(1+\th(\mu\xi)))},
\label{Q}
\end{equation}
where  $\l,\,\mu$ are positive constants defined by boundary
conditions as $\eta\to-\infty$ and $\xi\to-\infty$; $\rho$ is free
complex parameter.

This solution decreases with respect to spatial variables
exponentially. Besides if $\s=1$ and ${\mu\l\over4}|\rho|^2>1$ this
solution has singularities at some  lines $\eta=\const$ and
$\xi=\const$.

The solution (\ref{Q}) was called dromion in work \cite{F-S}. The
inverse scattering method for dromion-like solution of (\ref{ds1})
was developed in \cite{F-S}. From the results  of \cite{GK1},
\cite{GK2} and \cite{OK1} it follows that the soliton of the DS-2
equations is unstable with respect to small perturbation of the
initial data. These results stimulate the studies of the dromion
perturbation.

To construct the asymptotic solution over the large times we suppose
the parameter $\rho$ depending on a slow time  variable $\t=\ve t$
and obtain an explicit formula for the $\rho(\t)$.

The possibility of the full investigation of the linearized DS-1
equations plays the main role to construct the perturbation theory
of the nonlinear DS-1 equations. For (1+1)-models it has been found
in \cite{K1}, that spectral functions of the Lax pair form the basis
of solutions of the linearized nonlinear PDEs. The same  idea is
true for the DS-1 and DS-2 equations, as shown in \cite{OK5},
\cite{OK4}  and it is used below.

The  asymptotic analysis given in this paper is valid for the
solutions (\ref{Q}) without the singularities. It means if  $\s=1$,
then ${\mu\l\over4}|\rho(\t)|^2<1$. If the coefficient of the
perturbation $A>0$, then $|\rho|$ increases with respect to slow
time. It allows to say that the singularity may appear in the
leading term  of the asymptotics as $\t\to\infty$. However we can't
say this strongly for perturbed dromion in our situation, because
our asymptotics is usable only when $\t\ll\log(\log(\ve^{-1}))$. In
general case the appearance of the singularities in the solution of
non\-in\-teg\-rab\-le cases of the Davey-Stewartson equations is
known phenomenon \cite{PSSW}.

The obtained result for the problem about the interaction between
the long and short waves on the liquid surface shows when the depth
decreases the formal asymptotic solution can be described by
adiabatic perturbation theory of the dromion for the focusing  DS-1
equations at least for $|t|\ll\ve^{-1}\log(\log(\ve^{-1}))$.

The contents of the various sections are as follows. Section 2
contains a statement of a problem and a result. Section 3 is a brief
treatment of solving of linearized DS-1 equations using the basis
formed by solutions of the Dirac equation. In Section 4 we apply the
results of Section 3 and construct the first  correction of an
asymptotic  expansion for the solution of the problem formulated in
Section 1. Section 5 is devoted to reducing of a modulation equation
for the dromion parameter. In appendix we demonstrate explicit
formulas for functions which  are used to obtain the basis set
solving of the linearized DS-1 equations.

\section{Problem and result}

We construct the asymptotic solution of equ\-a\-ti\-ons (\ref{ds1})
on $\mod(O(\ve^2))$ uni\-form\-ly over large $t$. The perturbation
operator is $F\equiv AQ$ and the boundary conditions for $G_1$ and
$G_2$ are:
\begin{equation}
G_1|_{\xi\to-\infty}=u_1\equiv{\l^2\over2\cosh^2(\l\eta)},\quad
G_2|_{\eta\to-\infty}=u_2\equiv{\mu^2\over2\cosh^2(\mu\xi)}.
\label{bc}
\end{equation}
We seek the asymptotic solution as a sum of two first terms of
asymptotic expansions:
\begin{gather}
Q(\xi,\eta,t,\ve)=W(\xi,\eta,t,\t)+\ve U(\xi,\eta,t,\t),
\notag\\
G_1(\xi,\eta,t,\ve)=g_1(\xi,\eta,t,\t)+\ve V_1(\xi,\eta,t,\t),
\label{as1}\\
G_2(\xi,\eta,t,\ve)=g_2(\xi,\eta,t,\t)+\ve V_2(\xi,\eta,t,\t),
\notag
\end{gather}
where $\t=\ve t$ is slow time. The leading term of the
asymptotics has the form:
\[
W(\xi,\eta,t,\t)=q(\xi,\eta,t;\rho(\t)),
\]
and $g_1,\,\,g_2$ are:
\begin{gather*}
g_1(\xi,\eta,t,\t)=u_1-{\s\over2}\int_{-\infty}^{\xi}d\xi'\pt_{\eta}
|W(\xi',\eta,t,\t)|^2,
\\
g_2(\xi,\eta,t,\t)=u_2-{\s\over2}\int_{-\infty}^{\eta}d\eta'\pt_{\xi}
|W(\xi,\eta',t,\t)|^2.
\end{gather*}

Denote by $\g(\t)=1-\s{\mu\l\over4}|\rho(\t)|^2$  and $\g_0=\g(0)$.
The final result is given by
\begin{theorem}
If
\[
\g(\t)=\g_0^{\exp(2A\t)},\quad Arg(\rho(\t))\equiv\const,
\]
where $\g_0>1$ at $\s=-1$ and $0<\g_0<1$ at $\s=1$, then the
asymptotic solution (\ref{as1}) with respect to $\mod(O(\ve^2))$ is
useful uniformly over $t=O(\ve^{-1})$.
\end{theorem}
\begin{remark}
When the time is larger than $\ve^{-1}$, namely,
$t\ll \ve^{-1}\log(\log(\ve^{-1}))$, the formulas (\ref{as1}) are
asymptotic solution of (\ref{ds1}) with  respect to $\mod(o(1))$
only.
\end{remark}
\begin{remark}
D. Pelinovsky notes to author, that the modulation
of the parameter $|\rho|$ may  be obtained out of the ``energetic
equality'' \cite{DP}:
\[
\pt_t\int\int_{\Real^2}d\xi d\eta|Q|^2=
\ve\int\int_{\Real^2}d\xi d\eta(Q\bar F-\bar Q F).
\]
\end{remark}

\section{Solution of linearized equations}

In the present  section we obtain the formulas for solution of the
linearized DS-1 equations on the dromion  as a background:
\begin{equation}
\begin{split}
&i\pt_t U+(\pt_\xi^2+\pt_\eta^2)U+(G_1+G_2)U+(V_1+V_2)Q=iF
\\
&\pt_\xi V_1=-{\s\over2}\pt_\eta(Q\bar U+\bar Q U),\quad
\pt_{\eta}V_2={-\s\over2}\pt_\xi(Q\bar U+\bar Q U).
\end{split}
\label{l-ds}
\end{equation}

The results of inverse scattering transform \cite{F-S} and the set
of the basic functions \cite{OK5} are used for solving of the
linearized equations by Fourier method. It is reminded in Subsection
3.1.

However in contrast to the work \cite{OK5} here the solution of the
DS-1 equations with nonzero boundary conditions is considered. It
leads to changing of the dependency of the scattering data with
respect to time (see also \cite{F-S}) and of the formulas which
define the dependency  of Fourier coefficients of the solution for
the linearized DS-1 equations in contrast to obtained in \cite{OK5}.
This is explained in Subsection 3.2.

\subsection{Basic set for solving of the linearized DS-1 equation}

In the inverse scattering transform one use the matrix solution of
the Dirac  system to solve the DS-1 equations (see
\cite{Nizhnik}-\cite{F-S}):
\begin{equation}
\Bigg(
\begin{array}{cc}
\pt_\xi & 0\\
0 & \pt_\eta
\end{array}\Bigg)
\psi=-{1\over2}\Bigg(
\begin{array}{cc}
0 & Q\\ \s\bar Q & 0
\end{array}\Bigg)\psi.
\label{de}
\end{equation}

Let $\psi^+$ and $\psi^-$ be the matrix solutions of the Goursat
problem for the Dirac system with the boundary conditions (see
\cite{F-S}):
\begin{equation}
\begin{array}{ll}
\psi^+_{11}|_{\xi\to-\infty}=\exp(ik\eta), & \psi^+_{12}|_{\xi\to-\infty}=0,\\[1ex]
\psi^+_{21}|_{\eta\to\infty}=0,  & \psi^+_{22}|_{\eta\to-\infty}=\exp(-ik\xi);\\[1ex]
\psi^-_{11}|_{\xi\to-\infty}=\exp(ik\eta), & \psi^-_{12}|_{\xi\to\infty}=0,\\[1ex]
\psi^-_{21}|_{\eta\to-\infty}=0, & \psi^-_{22}|_{\eta\to-\infty}=\exp(-ik\xi).
\end{array}
\label{b-psi}
\end{equation}

Denote by $\psi^+_{(j)},\, j=1,2,$ the columns  of the matrix
$\psi^+$. The column $\psi_{(1)}^+$ is the solution of two systems
of equations. An additional system of time evolution is
\begin{gather}
\pt_t \psi^+_{(1)}=ik^2\psi^+_{(1)}+
i\Big(\begin{array}{cc} 1 & 0\\0 & -1\end{array}\Big)
(\pt_\xi-\pt_\eta)^2\psi^+_{(1)}
\notag
\\
\qquad{}+ i\Big(\begin{array}{cc} 0 & Q\\\s\bar Q  &
0\end{array}\Big) (\pt_\xi-\pt_\eta)\psi^+_{(1)}+
\Big(
\begin{array}{cc} iG_1 & -i\pt_\eta Q\\i\s\pt_\xi \bar Q & -iG_2\end{array}
\Big)
\psi^+_{(1)}.
\label{dtpsi}
\end{gather}
One can obtain the equation like this for  the  other columns of the
matrices  $\psi^{\pm}$. These equations are differ from
(\ref{dtpsi}) by the sign of $ik^2$ in first term of the right hand
side only.

Below we write two bilinear  forms defining analogues of the direct and
inverse Fourier transforms. First bilinear form is
\begin{equation}
(\chi,\mu)_f=\int_{-\infty}^{\infty}\int_{-\infty}^{\infty}d\xi
d\eta (\chi_1\mu_1\,\s\bar f\,+\,\chi_2\mu_2\,f).
\label{f1}
\end{equation}
Here $\chi_i$ and $\mu_i$ are the elements of the columns
$\chi$ and $\mu$.

Denote by $\phi_{(1)}$ and $\phi_{(2)}$ the solutions
conjugated to $\psi^+_{(1)}$ and $\psi^-_{(2)}$ with
respect to the bilinear form  (\ref{f1}).

Using  the formulas for the scattering  data (\cite{F-S}) one can
write these data as:
\begin{gather}
s_1(k,l)={1\over4\pi}(\psi^+_{(1)}(\xi,\eta,k),E_{(1)}(il\xi))_Q,
\label{sd1}
\\
s_2(k,l)={1\over4\pi}(\psi^-_{(2)}(\xi,\eta,k),E_{(2)}(il\eta))_Q.
\label{sd2}
\end{gather}
Here $E(z)={\hbox{diag}}(\exp(z),\exp(-z))$.

It is shown (\cite{F-S}) the elements of the matrices $\psi^{\pm}$
are analytic functions with respect to the variable $k$ when $\pm
Im(k)>0$. Using the scattering data one can write the nonlocal
Riemann-Hilbert problem for the $\psi^-_{11}$ and $\psi^+_{12}$ on
the real axes (\cite{F-S}):
\begin{gather*}
\psi^-_{11}(\xi,\eta,k)=\exp(ik\eta)+\exp(ik\eta)\bigg(\exp(-ik\eta)
\int_{-\infty}^{\infty} dl\,
s_1(k,l)\psi^+_{12}(\xi,\eta,l)\bigg)^-,
\\
\psi^+_{12}(\xi,\eta,k)=\exp(-ik\xi)\bigg(\exp(ik\xi)
\int_{-\infty}^{\infty} dl\,
s_2(k,l)\psi^-_{11}(\xi,\eta,l)\bigg)^+.
\end{gather*}
Here
\[
\bigg(f(k)\bigg)^{\pm}={1\over2i\pi}\int_{-\infty}^{\infty}
{dk'\,f(k')\over k'-(k\pm i0)}.
\]

The Riemann-Hilbert problem  for $\psi^-_{21}$ and
$\psi^+_{22}$ has the form:
\begin{gather*}
\psi^-_{21}(\xi,\eta,k)=
\exp(ik\eta)\bigg(\exp(-ik\eta)
\int_{-\infty}^{\infty} dl\,
s_1(k,l)\psi^+_{22}(\xi,\eta,l)\bigg)^-,
\\
\psi^+_{22}(\xi,\eta,k)=
\exp(-ik\xi)
+\exp(-ik\xi)\bigg(\exp(ik\xi)
\int_{-\infty}^{\infty} dl\,
s_2(k,l)\psi^-_{21}(\xi,\eta,l)\bigg)^+.
\end{gather*}

Introduce second bilinear form
\begin{equation}
\langle\chi,\mu\rangle_s=\int_{-\infty}^{\infty}\int_{-\infty}^{\infty}dk dl
(\chi^1(l)\mu^1(k)\,s_2(k,l)\,+\,\chi^2(l)\mu^2(k)\,s_1(k,l)),
\label{f2}
\end{equation}
where $\chi^{j}$ is the element of the row $\chi$.

Denote by $\varphi^{(j)},\, j=1,2,$  the row conjugated to
$\psi^{(j)}=[\psi^-_{j1},\psi^+_{j2}]$ with respect to the bilinear  form
(\ref{f2}). Formulate the result about the de\-com\-po\-si\-ti\-on obtained in
\cite{OK5}.

\begin{theorem}
Let $Q$ be such that
 $\pt^\a Q\in L_1\cap C^1$ for $|\a|\le
3$, if a function  $f$ is such that $\pt^\a f(\xi,\eta)\in L_1\cap
C^1$ for $|\a|\le4$, then  one can  represent  the $f$ in the form
\begin{equation}
f={-1\over\pi}\langle\psi^{(1)}(\xi,\eta,l),
\varphi^{(1)}(\xi,\eta,k)\rangle_{\hat f}, \label{f}
\end{equation}
where
\begin{equation}
\hat
f={1\over4\pi}(\psi^+_{(1)}(\xi,\eta,k),\phi_{(1)}(\xi,\eta,l))_f.
\label{hat-f}
\end{equation}
\end{theorem}

\subsection{Solving the linearized equation using Fourier method}

In  the preceding subsection we have state that the bilinear  forms
(\ref{f1}) and (\ref{f2}) may be used to decomposition like as a
Fourier integrals. The present subsection shows how to use this
decomposition for solving  of linearized DS-1 equation.

\begin{theorem}
Let $Q$ be the solution of the DS-1
equations with the boundary  conditions
$G_1|_{\xi\to-\infty}=u_1$ and
$G_2|_{\eta\to-\infty}=u_2,$ and $Q$ satisfies the
conditions of Theorem  2, the solution of the first of the
linearized DS-1 equation is smooth and integrable function
$U$ with respect to $\xi$ and $\eta$, where $\pt^\a U\in
L_1\cap C^1$ and $\pt^\a F\in L_1\cap C^1,$ for $|\a|\le
4$ and $t\in[0,T_0]$. Then
\begin{equation}
\pt_t\hat U=i(k^2+l^2)\hat U+\int_{-\infty}^{\infty} dk'
\hat U(k-k',l,t)\chi(k')+
\int_{-\infty}^{\infty} dl'\hat U(k,l-l',t)\kappa (l')+\hat F.
\label{u1}
\end{equation}
\end{theorem}

If the boundary conditions for the solution of the DS-1 equation
equal to zero, then $\chi\equiv\kappa\equiv0$. In this case the
formulas of the Theorem 3 allow to solve the linearized DS-1
equation in the explicit form. It was be done in \cite{OK5}. In
contrast of \cite{OK5}, we consider here the solution of the DS-1
equation with nonzero boundary conditions. It leads to integral
terms in the formula (\ref{u1}). In order to solve the linearized
DS-1 equation we must transform the  formula (\ref{u1}). In the
right hand side of (\ref{u1}) the integral terms are the
convolutions. Go over to the equations for the Fourier transform of
the functions $\hat U(k,l,t,\t)$ with respect to variables $k$ and
$l$. As a result we obtain  the linear Schr\"odinger equation:
\begin{equation}
i\pt_t \tilde U+(\pt_\xi^2+\pt_\eta^2)\tilde U+
(u_2(\xi\mu)+u_1(\eta\l))\tilde U=\tilde F.
\label{lsh}
\end{equation}
Here
\begin{equation}
\begin{split}
\tilde U(\xi,\eta,t,\t)={1\over2\pi}\int_{\Real^2}dkdl
\hat U(k,l,t,\t)\exp(-ik\eta-il\xi);
\\
\tilde F(\xi,\eta,t,\t)={1\over2\pi}\int_{\Real^2}dkdl
\hat F(k,l,t,\t)\exp(-ik\eta-il\xi).
\end{split}
\label{tilde-f}
\end{equation}
The same equation without the right hand side (as $\tilde F\equiv0$)
was obtained in \cite{F-S} for the time evolution of the scattering
data for the DS-1 equation.

One can construct the solution of the Cauchy problem for the
equations (\ref{lsh}) with the initial condition $\tilde
U(\xi,\eta)=0$ by the separation of the variables. In our case the
solution of the equations (\ref{lsh}) obtained by the Fourier method
has the form:
\begin{gather}
\tilde U(\xi,\eta,t)=
{1\over2\pi}\int_{\Real^2}dmdn\breve U(m,n,t)X(\xi,m)Y(\eta,n)
\exp(-it(m^2+n^2))\notag\\
\qquad{}+{1\over\sqrt{2\pi}}
\int_{\Real}dn\breve U_{\mu}(n,t)Y(\eta,n)X_{\mu}(\xi)\exp(-it(n^2-\mu^2))
\notag\\
\qquad{}+{1\over\sqrt{2\pi}}
\int_{\Real}dm\breve U_{\l}(m,t)X(\xi,m)Y_{\l}(\eta)\exp(-it(m^2-\l^2))
\notag\\
\qquad{}+\breve U_{\mu,\l}X_{\mu}(\xi)Y_{\l}(\eta)\exp(it(\mu^2+\l^2)).
\label{tilde-U}
\end{gather}
Here we use notations:
\begin{gather}
X(m,\xi)={\mu\th(\mu \xi)+im\over im-\mu}\exp(-im\xi),
\quad X_{\mu}(\xi)={1\over2\ch(\mu\xi)};
\notag\\
Y(n,\eta)={\l\th(\l \eta)+in\over in-\l}\exp(-in\eta),
\quad Y_{\l}(\eta)={1\over2\ch(\l\eta)};
\notag
\\[1ex]
\begin{split}
&\pt_t\breve U=i(m^2+n^2)\breve U+\breve F(m,n,t),
\quad
\pt_t \breve U_{\mu}=i(n^2-\mu^2)\breve U+\breve F_{\mu}(n,t),\\
&\pt_t \breve U_{\l}=i(m^2-\l^2)\breve U+\breve F_{\l}(m,t),
\quad
\pt_t \breve U_{\mu\l}=-i(\l^2+\mu^2)\breve U+\breve F_{\mu\l}(t),\\
&\breve U|_{t=0}=\breve U_{\mu}|_{t=0}=\breve U_{\l}|_{t=0}=
\breve U_{\mu\l}|_{t=0}=0;
\end{split}
\label{Cauch}
\\[1ex]
\begin{split}
&\breve F(m,n,t)=\int_{\Real^2}d\xi d\eta \tilde F(\xi,\eta,t)
\bar X(m,\xi)\bar Y(n,\eta),\\
&\breve F_{\mu}(n,t)=\int_{\Real^2}d\xi d\eta \tilde F(\xi,\eta,t)
\bar X_{\mu}(\xi)\bar Y(n,\eta),\\
&\breve F_{\l}(m,t,\t)=\int_{\Real^2}d\xi d\eta \tilde F(\xi,\eta,t)
\bar X(m,\xi)\bar Y_{\l}(\eta),\\
&\breve F_{\mu\l}(t,\t)=\int_{\Real^2}d\xi d\eta \tilde F(\xi,\eta,t)
\bar X_{\mu}(\xi)\bar Y_{\l}(\eta).
\end{split}
\label{breve-f}
\end{gather}

To solve the linearizing equations we must apply the formulas
(\ref{hat-f}), (\ref{tilde-f}), (\ref{breve-f}), solve the Cauchy
problems (\ref{Cauch}) and apply formulas (\ref{tilde-U}), inverse
Fourier transform to $\tilde U$ and then apply formula (\ref{f}) for
$\hat U$. As a result we obtain the solution of the linearized DS-1
equations. The functions $V_{1,2}$ may be obtained by direct
integrating of second and third equations from (\ref{l-ds}).

\section{Equation for the first correction}

This and next section contain a proof of the Theorem 1.  In this
part the equation for the slow modulation of the parameter
$\rho(\t)$ is obtained. This equation is necessary and sufficient
condition for the uniform boundedness of the  first correction of
the expansion (\ref{as1}) over $t=O(\ve^{-1})$.

Some complication appears when one use the Fourier method from
preceding section to solve the linearized DS-1 equations on dromion
as a background with slow varying parameter $\rho(\t)$. In this case
the basic functions used for solving the linearized equations depend
on slow variable $\t$ also and they are the basic set of asymptotic
solutions with respect to $\mod(O(\ve))$ only.

Substitute  the formula (\ref{as1}) into the equations (\ref{ds1}).
Equate the coefficients with the same power of $\ve$. The equations
as $\ve^0$ are realized since $W,g_1,g_2$ are the  asymptotic
solution of nonperturbed DS-1 equations. For the first correction we
obtain the linearized DS-1 equations:
\begin{equation}
\begin{split}
&i\pt_t
U+(\pt_\xi^2+\pt_\eta^2)U+(g_1+g_2)U+(V_1+V_2)W= iH
\\
&\pt_\xi
V_1=-{\s\over2}\pt_\eta(W\bar U+\bar W U),\quad
\pt_{\eta}V_2={-\s\over2}\pt_\xi(W\bar U+\bar W U),
\end{split}
\label{lds1}
\end{equation}
where
\[
H=AW-\pt_\t W.
\]

Before to use the formulas from the preceding section we reduce the
form of the right hand side in first of the equations (\ref{lds1}).
In the leading  term the parameter $\rho$ depends on the $\t$ only.
The other parameters of the dromion ($\l$ and $\mu$) depend only on
the boundary conditions and have not changes under perturbation of
the equations. For more convenience we represent
$\rho(\t)=r(\t)\exp(i\a(\t))$ where $r(\t)=|\rho(\t)|$ and
$\a(t)={\hbox{Arg}}\rho(\t)$. The derivation of $W$ with respect to
slow variable $\t$ can be written as:
\[
\pt_\t W=\pt_r W  r'+\pt_\a W \a'.
\]
Here the derivatives $r'$ and $\a'$ are unknown yet.

Compute the function $\hat H$. Using  the Theorem 2 and formulas for
the functions $\psi_{+}$ and $\phi$ (see Appendix) we obtain:
\[
\widehat H(k,l,t,\t)=\exp(-it(\l^2+\mu^2))\bigg(\widehat P(k,l;\rho)-
\widehat R(k,l;\rho)\pt_\t\rho\bigg),
\]
where
\[
\widehat
P(k,l;\rho)=\exp(it(\l^2+\mu^2))\widehat{AW},\quad \widehat R(
k,l;\rho)=\exp(it(\l^2+\mu^2))\widehat{\pt_\t W}.
\]
In these formulas we have write the dependence on time in explicit
form. It allows to remove the secular terms (increasing with respect
to $t$) from the asymptotic solution (\ref{as1}) using modulation of
the parameters $r(\t)$ and $\a(\t)$.

The differential equations for $\breve{U}$ have the forms:
\begin{gather*}
\pt_t\breve U=i(m^2+n^2)\breve U+
\exp(-it(\l^2+\mu^2))(\breve P(m,n;\rho)-\breve R(m,n;\rho)),
\\
\pt_t \breve U_{\mu}=i(n^2-\mu^2)\breve U+
\exp(-it(\l^2+\mu^2))(\breve P_\mu(n;\rho)-\breve R_\mu(n;\rho)),\\
\pt_t \breve U_{\l}=i(m^2-\l^2)\breve U+
\exp(-it(\l^2+\mu^2))(\breve P_{\l}(m;\rho)-\breve R_\l(m;\rho)),
\\
\pt_t \breve U_{\mu\l}=-i(\l^2+\mu^2)\breve U+
\exp(-it(\l^2+\mu^2))(\breve P_{\mu\l}(\rho)-\breve R_{\mu\l}(\rho)).
\end{gather*}

Solutions of these equations are
\begin{gather*}
\breve U(m,n,t)={\mu^2+\l^2\over n^2+m^2+\mu^2+\l^2}
\bigg(\breve P(m,n;\rho)-
\breve R(m,n;\rho)\bigg)\exp(-it(\mu^2+\l^2)),\\
\breve U_\l(m,t)={\mu^2+\l^2\over m^2+\mu^2}\bigg(\breve P_\l(m;\rho)-
\breve R_\l(m;\rho)\bigg)\exp(-it(\mu^2+\l^2)),\\
\breve U_\mu(n,t)={\mu^2+\l^2\over n^2+\l^2}\bigg(\breve P_\mu(n;\rho)-
\breve R_\mu(n;\rho)\bigg)\exp(-it(\mu^2+\l^2)),\\
\breve U_{\mu\l}(t)=\bigg(\breve P_{\mu\l}(m;\rho)- \breve
R_{\mu\l}(m;\rho)\bigg)t\exp(-it(\mu^2+\l^2)).
\end{gather*}

One can see that the secular terms may appear because of the last
term in the equation for $\breve U_{\mu\l}$. The eliminating  of this
term leads us to the equation for $\rho(\t)$:
\begin{equation}
\breve R_{\mu\l}-\breve P_{\mu\l}=0,\quad
\rho|_{\t=0}=\rho_0.
\label{sec}
\end{equation}

As a  result  $\breve U_{\mu\l}\equiv 0$. The another solutions
$\breve U(m,n,t),\,\breve U_\mu(n,t)$ and $\breve U_\l(m,t)$ are
bounded with  respect to all arguments and over all times.

One must return to the original of the images $\breve U$ to say
about boundedness of the solution $U, V_1, V_2$ for the equations
(\ref{lds1}). One can  see the direct (from $U$ into $\breve U$) and
inverse (from $\breve U$ into $U$) integral transforms as the
Fourier transform from the smooth  and exponentially decreasing
functions with respect to the corresponding variables. The Fourier
transform moves such functions into analytic functions near the real
axis. The inverse transform moves these analytic functions into the
exponential decreasing functions.  So the solution of (\ref{lds1})
is bounded and decreasing exponentially  with respect to the spatial
variables.

\section{Modulation equation}

Here the equation (\ref{sec}) for the parameter $\rho(\t)$
is reduced to the more  con\-vi\-ent form. Write the
derivative of the leading term with respect to the slow
variable $\t$.
\[
i\pt_\t W=-\a'W- iW{r'\over  r}+{2iW\over 1-\s
r^2{\mu\l\over16}(1+\th(\mu\xi))(1+\th(\l\eta))}\,{r'\over
r}.
\]
Denote by $\G={\mu\l\over4}r^2$ and compute the images
$\breve{(\cdot)}_{\mu\l}$ of every term.
\begin{gather*}
\breve{(W)}_{\mu\l}={\s\bar\rho\exp(-it(\l^2+\mu^2))\over 8}
(\s\G-1)\log|1-\s\G|;
\\
\breve{(iW)}_{\mu\l}=-i\s\G{\rho\exp(-it(\l^2+\mu^2))\over 8}.
\end{gather*}
Denote the image of the last term by $\breve h_{\mu\l}$.
Its image has the form:
\[
\breve h_{\mu\l}={r'\over r}{\s\bar\rho\exp(-it(\l^2+\mu^2))\over8}
(1-\s\G)\bigg({1\over 1-\s\G}-1-\log|1-\s\G|\bigg).
\]
The image $\breve{(\cdot)}$ of $AW$ has the similar form.

Substitute these formulas into ({\ref{sec}}) and separate the real and
the imaginary parts of this equation, then
\begin{gather*}
\a'=0,\\
{r'\over r}(\s\G-1)\log|1-\s\G|- {r'\over
r}\bigg(\s\G-(1-\s\G)\log|1-\s\G|\bigg)\\
\qquad{}+A(\s\G-1)\log|1-\s\G|=0.
\end{gather*}
Use  the notation for  $\G$, then the second equation has
the form:
\[
{d\G\over d\t}=-2\s A (1-\s\G)\log|1-\s\G|.
\]
This equation defines the evolution  of the absolute value
of the complex parameter $\rho$. The argument of this
parameter do not change under the per\-tur\-ba\-ti\-on
$F=AQ$.

Denote by $\g=1-\s\G$, rewrite the equation for $\g$. Then we obtain:
\[
\g'=2A\g\log(\g).
\]
The solution for this equation  has the form:
\[
\g(\t)=\exp(C\exp(2A\t)),
\]
where $\g|_{\t=0}=\exp(C)$, then we can write the $\g(\t)$
in  the form:
\[
\g(\t)=\g_0^{\exp(2A\t)}.
\]
The Theorem 1 is proved.

\subsection*{Acknowledgements}
I thank D. Pelinovsky
for the discussions of the results and for the helpful remarks. Also
I thank V.Yu Novokshenov, S. Glebov and N. Enikeev. Their remarks had
allowed to improve this paper. This work was supported by RFBR (00-01-00663,
00-15-96038) and INTAS (99-1068).

\appendix

\section*{Appendix}
\renewcommand{\thesection}{A}
\subsection{Explicit formulas}

Here  we remain the explicit forms of the solution for the
Dirac equation with the dromion-like potential. These
forms were obtained in \cite{F-S}. In  our computations we
use first column of the matrix $\psi^+$ only.
\begin{gather}
\left(\begin{array}{c}
\psi^+_{11}\\
\psi^+_{21}\end{array}\right)
=\left(\begin{array}{c}\exp(ik\eta)\\0\end{array}\right)+
\frac{\int^{\infty}_\eta {dp\l\exp(ikp))\over2\cosh(\l p)}}{
(1-\s|\rho|^2{\l\mu\over16}(1+\tanh(\mu\xi))(1+\tanh(\l\eta)))}
\notag\\
\qquad{}\times
\left(\begin{array}{c}
{-\s|\rho|^2\l\mu(1+\tanh(\mu\xi))\over8\cosh(\l\eta)}\\
{\s\bar\rho\mu\exp(-it(\l^2+\mu^2))\over2\cosh(\mu\xi)}
\end{array}\right).
\label{psi+1}
\end{gather}

\subsection{Solutions  of  conjugated equations}

Here we  write the problems for the functions conjugated
to $\psi^{\pm}$ with respect to the bilinear forms.

The matrix $\phi$ is the solution of the boundary problem
con\-ju\-ga\-ted to the solutions of the problem (\ref{de}),
(\ref{bc}) with respect to the bilinear form (\ref{f1}). First
column of the matrix $\phi$ is the solution of the integral equation:
\begin{gather*}
\phi_{11}(\xi,\eta,l,t)=\exp(il\xi)+
{1\over2}\int_{-\infty}^\eta d\eta'
Q(\xi,\eta',t)\phi_{21}(\xi,\eta',l,t),\\
\phi_{21}(\xi,\eta,l,t)=
{-1\over2}\int_{\xi}^{\infty} d\xi' \bar
Q(\xi',\eta,t)\phi_{11}(\xi',\eta,l,t).
\end{gather*}

The explicit formula for first column of the matrix $\phi$ used in
Section 4 for the dromion potential has the form:
\begin{gather}
\left(\begin{array}{c}
\phi_{11}\\
\phi_{21}\end{array}\right)
=\left(\begin{array}{c}\exp(il\xi)\\0\end{array}\right)+
{\int^{\infty}_\xi {dp\mu\exp(ilp))\over2\cosh(\mu
p)}\over
(1-\s|\rho|^2{\l\mu\over16}(1+\tanh(\mu\xi))(1+\tanh(\l\eta)))}
\notag\\
\qquad{}\times
\left(\begin{array}{c}
{-\s|\rho|^2\l\mu(1+\tanh(\l\eta))\over8\cosh(\mu\xi)}\\
{-\s\bar\rho\l\exp(-it(\l^2+\mu^2))\over2\cosh(\l\eta)}
\end{array}\right).
\label{phi1}
\end{gather}

Second bilinear  form (\ref{f2}) allows to  write the integral
equations conjugated to the integral equations which were obtained
from the nonlocal Riemann-Hilbert equation for the matrices
$\psi^{\pm}$ in \cite{F-S}. Using the equation for the matrix
$\varphi(\xi,\eta,l)$ one can show that  $\varphi_{11}$ and
$\varphi_{21}$ are the solutions of the boundary problem for the
Dirac system or the equivalent integral equations:
\begin{gather*}
\varphi_{11}(\xi,\eta,l)=\exp(il\xi)+
{1\over2}\int_{-\infty}^\eta d\eta'Q(\xi,\eta',t)\varphi_{21}(\xi,\eta',l),\\
\varphi_{21}(\xi,\eta,l)=
{1\over2}\int_{-\infty}^\xi d\xi'\bar Q(\xi',\eta,t)\varphi_{11}(\xi',\eta,l).
\end{gather*}

The functions $\varphi_{11}$ and  $\varphi_{21}$ have the similar
form as the functions $\phi_{11}$ and  $\phi_{21}$.

\label{lastpage}

\end{document}